\documentclass[a4paper,11pt,DIV=11,abstract=on]{scrartcl}
\usepackage{mathtools}
\usepackage{amsmath, amsthm, amssymb}
\usepackage{fontenc}
\usepackage[utf8]{inputenc}
\usepackage[english]{babel}
\usepackage{graphicx}
\usepackage{subcaption,xargs}
\usepackage{algorithm, booktabs}
\usepackage[noend]{algpseudocode}
\usepackage{enumitem,listings,microtype}
\usepackage{environ,ifthen,xcolor}
\usepackage{multirow}
\usepackage[
  backend=bibtex,
  style=numeric-comp,
  firstinits=true,
  natbib=true,
  maxcitenames=5,
  sortlocale=en_EN,
  url=true, 
  doi=true,
  eprint=true
]{biblatex}
\usepackage{hyperref}
\newcommand{\e}{{\,\mathrm{e}}}
\newcommand{\im}{{\mathrm{i}}}
\newcommand{\p}{{\mathrm{p}}}

\newcommand{\dx}{\,\mathrm{d}}
\makeatletter
\newcommandx{\abs}[2][1=\@empty]{#1\lvert #2 #1\rvert}
\newcommandx{\norm}[3][1=\@empty,3=\@empty]{#1\lVert #2 #1\rVert_{#3}}
\makeatother
\newcommand{\tT}{\mathrm{T}}

\newcommand{\vect}[1]{\mathbf{#1}}
\newcommand{\mat}[1]{\mathbf{#1}}
\newcommand{\Stiffness}{\mathcal{C}}
\newcommand{\eff}{\mathrm{eff}}

\newcommand{\N}{\mathbb{N}}
\newcommand{\Z}{\mathbb{Z}}
\newcommand{\R}{\mathbb{R}}

\newcommand{\T}{\mathbb{T}}

\newcommand{\tensorProd}{:}

\DeclareMathOperator{\Id}{Id}
\DeclareMathOperator{\Fourier}{\mathcal F}
\DeclareMathOperator{\generatingSet}{\mathcal G}
\DeclareMathOperator{\Grad}{\nabla}
\DeclareMathOperator{\Pattern}{\mathcal P}
\DeclareMathOperator{\GradSym}{\Grad}
\DeclareMathOperator{\spanOp}{span}

\DeclareMathOperator{\Translate}{\mathcal T}
\DeclareMathOperator{\diag}{diag}

\DeclareMathOperator{\Fundamental}{I}

\newtheorem{theorem}{Theorem}%[section]
\newtheorem{definition}[theorem]{Definition}%[section]

\title{Approximation of Periodic PDE Solutions with Anisotropic Translation
Invariant Spaces}

\author{%
Ronny Bergmann\thanks{%
Department of Mathematics,
University of Kaiserslautern,
Postfach 3049, D-67653 Kaiserslautern, Germany\newline$\{$bergmann, dmerkert$\}$@mathematik.uni-kl.de.
}
\and
Dennis Merkert\footnotemark[1]}
\date{May 02, 2017}

\begin{document}
\maketitle
\begin{abstract}
  \noindent\small
  We approximate the quasi-static equation of linear elasticity in translation
  invariant spaces on the torus. This unifies different FFT-based discretisation
  methods into a common framework and extends them to anisotropic lattices. We
  analyse the connection between the discrete solution spaces and demonstrate
  the numerical benefits. Finite element methods arise as a special case of
  periodised Box spline translates.
\end{abstract}

\section{Introduction}

The simulation of composite materials is the foundation of many modern tools and
applications. Typical simulated materials are fibre reinforced plastics or metal
foams where one is interested in the characterization of their macroscopic
elastic properties. These properties are described by the quasi-static equation of
linear elasticity in homogenization with periodic boundary conditions. For
this equation Moulinec and
Suquet~\cite{MoulinecSuquet1994,MoulinecSuquet1998} developed a FFT-based
algorithm which approximates the solution using truncated Fourier
series~\cite{Vondrejc2014}.

This manuscript uses anisotropic spaces of periodic
translates~\cite{LangemannPrestin2010WaveletAnalysis,Bergmann2013FFT} for the
approximation. It unifies and generalizes the approaches of Vond\v{r}ejc,
et.al~\cite{Vondrejc2014} and Brisard and
Dormieux~\cite{BrisardDormieux2010Framework}.
The detailed analysis can be found in~\cite{BergmannMerkert2017}.

The continuous equation allows for two equivalent formulations. They have the
same discrete solution in a space of translates if and only if the space yields
the truncated Fourier series approach which corresponds to using spaces of
translates generated by a Dirichlet kernel on anisotropic
patterns~\cite{BergmannMerkert2016}.

The remainder is structured as follows: We introduce patterns and anisotropic
spaces of translates and the partial differential equation in
Section~\ref{sec:preliminaries}. The approximation of the PDE and connections
between solution spaces stemming from the two formulations of the equation are
explained in Section~\ref{sec:approximation}. The effects of this
framework is then demonstrated numerically in Section~\ref{sec:numerics} followed
by a conclusion in Section~\ref{sec:conclusion}.

\section{Preliminaries}\label{sec:preliminaries}
We first introduce the spaces of interest, namely translation invariant spaces
of multivariate periodic square integrable functions \(f\colon\T^d\to\mathbb
C\), where \(\T^d\) denotes the \(d\)-dimensional torus~$\T^d \cong
[-\pi,\pi)^d$. The translation invariance is defined with respect to patterns
which generalize the usual rectangular sampling grids. On these patterns, a
discrete Fourier transform is defined which can be employed to characterize the
spaces of interest using Fourier coefficients of the involved functions.
Next, two equivalent formulations of the equation of quasi-static elasticity
with periodic boundary conditions are introduced.

\subsection{Patterns and translation invariant spaces and the Fourier transform}

For any regular matrix $\mat{M} \in \mathbb Z^{d\times d}$, we define the
congruence relation for $\vect{h},\vect{k} \in \mathbb Z^d$ with respect
to~$\mat{M}$ by
\begin{equation*}
  \vect{h} \equiv \vect{k} \bmod \mat{M}
  \Leftrightarrow \exists\,\vect{z} \in \mathbb Z^d\colon \vect{k} = 
  \vect{h} + \mat{M}\vect{z}\text{.}
\end{equation*}
We define the lattice
\[
  \Lambda(\mat{M}) \coloneqq \mat{M}^{-1}\mathbb Z^d
  = \{\vect{y}\in\mathbb R^d : \mat{M}\vect{y} \in \mathbb Z^d\},
\]
and the pattern \(\Pattern(\mat{M})\) as any set of congruence representants of
the lattice with respect to \(\bmod\ 1\), e.g.\ $ \Lambda(\mat{M})\cap[0,1)^d$
or~$\Lambda(\mat{M})\cap\bigl[-\tfrac{1}{2},\tfrac{1}{2}\bigr)^d$. For the rest
of the paper we will refer to the set of congruence class representants in the 
symmetric unit cube \(\bigl[-\tfrac{1}{2},\tfrac{1}{2}\bigr)^d\),
i.e.~$\Pattern(\mat{M}) = \Lambda(\mat{M})\cap
\bigl[-\tfrac{1}{2},\tfrac{1}{2}\bigr)$. The generating
set \(\generatingSet(\mat{M})\) is defined by \(\generatingSet(\mat{M}) 
\coloneqq \mat{M}\Pattern(\mat{M})\) for any pattern \(\Pattern(\mat{M})\). For
both, the number of elements is given by \(
  \abs{\Pattern(\mat{M})}
  =\abs{\generatingSet(\mat{M})}
  =\abs{\det{\mat{M}}}
  \eqqcolon m
\), cf.~\cite[Lemma II.7]{deBoorHoelligRiemenscheider1993BoxSplines}.

A space of square integrable functions \(V\subset L^2(\T^d)\) is called
\(\mat{M}\)-invariant, if for all \(\vect{y}\in\Pattern(\mat{M})\) and all
functions \(f\in V\) the translates \(\Translate({\vect{y}})f \coloneqq
f(\cdot-2\pi\vect{y})\in V\), i.e. especially
\[
  V_{\mat{M}}^{f} \coloneqq \spanOp\bigl\{\Translate(\vect{y})f\,;\,
  \vect{y}\in\Pattern(\mat{M})\bigr\}
\]
is \( \mat{M} \)-invariant.
A function \(g\in V_{\mat{M}}^f \) is of the form
\begin{equation}\label{eq:functionintranslates}
  g=\displaystyle\sum_{\vect{y}\in\Pattern(\mat{M})}
  a_{\vect{y}}\Translate(\vect{y})f,
\end{equation}
where \(\vect{a} = (a_{\vect{y}})_{y\in\Pattern(\mat{M})}\in\mathbb C^m\).

A function $\Fundamental_\mat{M} \in V^f_\mat{M}$ is called fundamental
interpolant if $\Fundamental_\mat{M}(2 \pi \vect{y}) = 1$ if $\vect{y} \equiv
\vect{0} \mod \mat{M}$ and $0$ else.

The discrete Fourier matrix on the pattern $\Pattern(\mat{M})$ is
defined~\cite{ChuiLi:1994} by
\begin{equation}\label{eq:Fouriermatrix}
  \mathcal F(\mat{M})
  \coloneqq
  \frac{1}{\sqrt{m}}
  \Bigl(
    \e^{- 2\pi \im \vect{h}^\tT\vect{y}}
  \Bigr)_{%
    \vect{h} \in \generatingSet(\mat{M}^\tT),\,%
    \vect{y} \in \Pattern(\mat{M})%
  },
\end{equation}
where $\vect{h}\in\generatingSet(\mat{M}^\tT)$ the rows and $\vect{y}
\in \Pattern(\mat{M})$ indicate the columns of the Fourier matrix~\(\mathcal
F(\mat{M})\), and the discrete Fourier transform is defined by
\begin{equation}\label{eq:FourierTransform}
  \vect{\hat a} = (\hat a_{\vect{h}})_{\vect{h}\in\generatingSet(\mat{M}^\tT)}
  = \mathcal F(\mat{M})\vect{a},
\end{equation}
where~$\vect{a} = (a_{\vect{y}})_{\vect{y}\in\Pattern(\mat{M})}\in\mathbb C^m$
and \(\vect{\hat a}\) are arranged in the same ordering as the columns and
rows in~\eqref{eq:Fouriermatrix}, respectively. A fast Fourier transform
can be derived, see~\cite[Theorem~2]{Bergmann2013FFT}.

For \( f\in L^2(\T^d) \) an easy calculation on the Fourier coefficients
\(c_{\vect{k}}(f) = \langle f,\e^{\im\vect{k}^\tT\circ}\rangle\),
\(\vect{k}\in\mathbb Z^d\), using the unique decomposition of \(
\vect{k}\in\mathbb Z^d \) into \( \vect{k} =
\vect{h}+\mat{M}^\tT\vect{z} \), \(
\vect{h}\in\generatingSet(\mat{M}^\tT),\ \vect{z}\in\mathbb Z^d \),
yields, that \( g\in V_{\mat{M}}^f \) holds if and only
if~\cite[Theorem 3.3]{LangemannPrestin2010WaveletAnalysis}
\begin{equation}\label{eq:inTranslates:ck}
  c_{\vect{h}+\mat{M}^\tT\vect{z}}(g) 
  = \hat a_{\vect{h}}c_{\vect{h}+\mat{M}^\tT\vect{z}}(f)
  \quad\text{for all }
  \vect{h}\in\generatingSet(\mat{M}^\tT),\ \vect{z}\in\mathbb Z^d\text{,}
\end{equation}
holds, where \(\vect{\hat a} = \mathcal F(\mat{M})\vect{a}\) for
\(\vect{a}\) from~\eqref{eq:functionintranslates}. Finally we define for
a generalized sequence \(\vect{a} = \{a_{\vect{k}}\}_{\vect{k}\in\mathbb
Z^d}\) and a regular integer matrix \(\mat{M}\in\mathbb Z^{d\times d}\)
the bracket sum
\begin{equation}\label{eq:bracketsum}
  \bigl[\vect{a}\bigr]_{\vect{k}}^{\mat{M}}
  \coloneqq
  \sum_{\vect{z}\in\mathbb Z^d} a_{\vect{k}+\mat{M}^\tT\vect{z}},
  \qquad \vect{k}\in\mathbb Z^d.
\end{equation}
The discrete Fourier coefficients \(c_{\vect{h}}^{\mat{M}}\) are related to
the Fourier coefficients for a function~\( f\in A(\T^d) \), where
	\(A(\T^d)\) denotes the Wiener Algebra, i.e. the space of functions
	with an absolutely convergent Fourier series. 

Special cases for the choice of the space of translated $V_\mat{M}^f$ are given
by the Dirichlet kernel $f = D_\mat{M}$ with
\begin{equation}
  c_\vect{k}(D_\mat{M}) = \begin{cases} 1, & \vect{k} \in \generatingSet(\mat{M}^\tT)\\
  0,& \text{else} \end{cases},\ \vect{k} \in \Z^d
\end{equation}
and the de la Vall\'ee Poussin means $f = f_{\mat{M},\alpha}$ given by
\begin{equation}
  c_\vect{k}(f_{\mat{M},\alpha}) = \frac{1}{\sqrt{m}} B_\Xi(\mat{M}^{-\tT}
  \vect{k}),\ \vect{k} \in \Z^d.
\end{equation}
The function $B_\Xi$ denotes a Box spline with $\Xi = (\diag(\alpha)\Id_d) \in
\R^{d \times 2d}$ and $\alpha \in [0,1]^d$.

\subsection{Periodic homogenization in linear elasticity}

FFT-based methods to solve the equations of periodic homogenization in linear
elasticity were first introduced by Moulinec and
Suquet~\cite{MoulinecSuquet1994,MoulinecSuquet1998} and later interpreted as a
Galerkin projection onto trigonometric sums by Vond\v{r}ejc
et.al.~\cite{Vondrejc2014}. The partial differential equation (PDE) is
formulated as a variational equation in terms of the symmetric strain
$\varepsilon \in L^2(\T^d)^{d \times d}$ with zero mean given a symmetric
macroscopic strain $\varepsilon^0 \in \R^{d \times d}$ and a stiffness
distribution $\Stiffness \in L^\infty(\T^d)^{d \times d \times d \times d}$.
This stiffness is assumed to be uniformly elliptic, i.e. there exists a constant
$c \in \R$ such that $\langle \Stiffness \tensorProd \gamma, \gamma \rangle
\geq c \norm{\gamma}^2$ for all $\gamma \in L^2(\T^d)^{d \times d}$ symmetric.
Further, assume that $\Stiffness$ has the symmetries
$\Stiffness_{ijkl} = \Stiffness_{jikl} = \Stiffness_{ijlk} = \Stiffness_{klij}$
for $i,j,k,l = 1,\dots,d$. 

The equation makes use of the Green operator $\Gamma^0$ defined by its action in
frequency domain
\begin{align}
  \label{eq:gamma}
    \Gamma^0 \tensorProd \varepsilon &\coloneqq
    \sum_{\vect{k} \in \Z^d}
    \hat{\Gamma}^0_\vect{k} \tensorProd c_\vect{k}(\varepsilon)
    e^{2 \pi \im \vect{k}^\tT \cdot}
    \intertext{and Fourier coefficients}
    \hat{\Gamma}^0_\vect{k} \tensorProd c_\vect{k}(\varepsilon)
    &\coloneqq
    \GradSym_\vect{k}
    \Bigl(
      \overline{\GradSym_\vect{k}}^\tT \tensorProd \Stiffness^0 \tensorProd
      \GradSym_\vect{k}
    \Bigr)^{-1}
    \overline{\GradSym_\vect{k}}^\tT c_\vect{k}(\varepsilon),\ \vect{k} \in
    \Z^d
    \intertext{and for $u \in L^2(\T^d)^d$}
    \GradSym_\vect{k} c_\vect{k}(u) &\coloneqq \frac{\im}{2} \bigl( \vect{k}
    c_\vect{k}(u)^\tT + c_\vect{k}(u) \vect{k}^\tT\bigr),\ \vect{k} \in \Z^d.
\end{align}
The positive definite reference stiffness $\Stiffness^0 \in \R^{d \times d
\times d \times d}$ with the same symmetries as $\Stiffness$ is a free
parameter and $\Gamma^0 \tensorProd \varepsilon$ denotes the product of the
fourth-order tensor $\Gamma^0$ with the second-order tensor $\varepsilon$.

With these definitions at hand there are two equivalent formulations of the PDE
which yield the strain $\varepsilon$, cf.~\cite{Vondrejc2014}:
The variational equation (VE) states
\begin{equation}
  \bigl\langle
  \Stiffness^0 \Gamma^0 \Stiffness \tensorProd \varepsilon,\gamma
  \bigr\rangle =
  -\bigl\langle
  \Stiffness^0 \Gamma^0 \Stiffness \tensorProd \varepsilon^0,\gamma
  \bigr\rangle
  \label{eq:weak_pde_projected}
\end{equation}
for all symmetric $\gamma \in L^2(\T^d)^{d \times d}$ and the Lippmann-Schwinger
equation (LS) requires
\begin{equation}
  \label{eq:lippmann-schwinger}
  \bigl\langle
  \varepsilon +
  \Gamma^0 \bigl(
    \Stiffness - \Stiffness^0
  \bigr) \tensorProd
  \bigl(
    \varepsilon + \varepsilon^0
  \bigr),
  \gamma
  \bigr\rangle = 0
\end{equation}
to hold true for all symmetric $\gamma \in L^2(\T^d)^{d \times d}$.

The effective stiffness matrix $\Stiffness^\eff$ is given by
\begin{equation}
  \Stiffness^\eff \tensorProd \varepsilon^0
  \coloneqq \int_{\T^d} \Stiffness \tensorProd \varepsilon \dx x
\end{equation}
where $\varepsilon$ solves the PDE~\eqref{eq:lippmann-schwinger} with
macroscopic strain $\varepsilon^0$.

\section{Approximation of periodic PDE solutions}\label{sec:approximation}

\begin{figure*}[t]\centering
  \includegraphics[width=.675\textwidth]{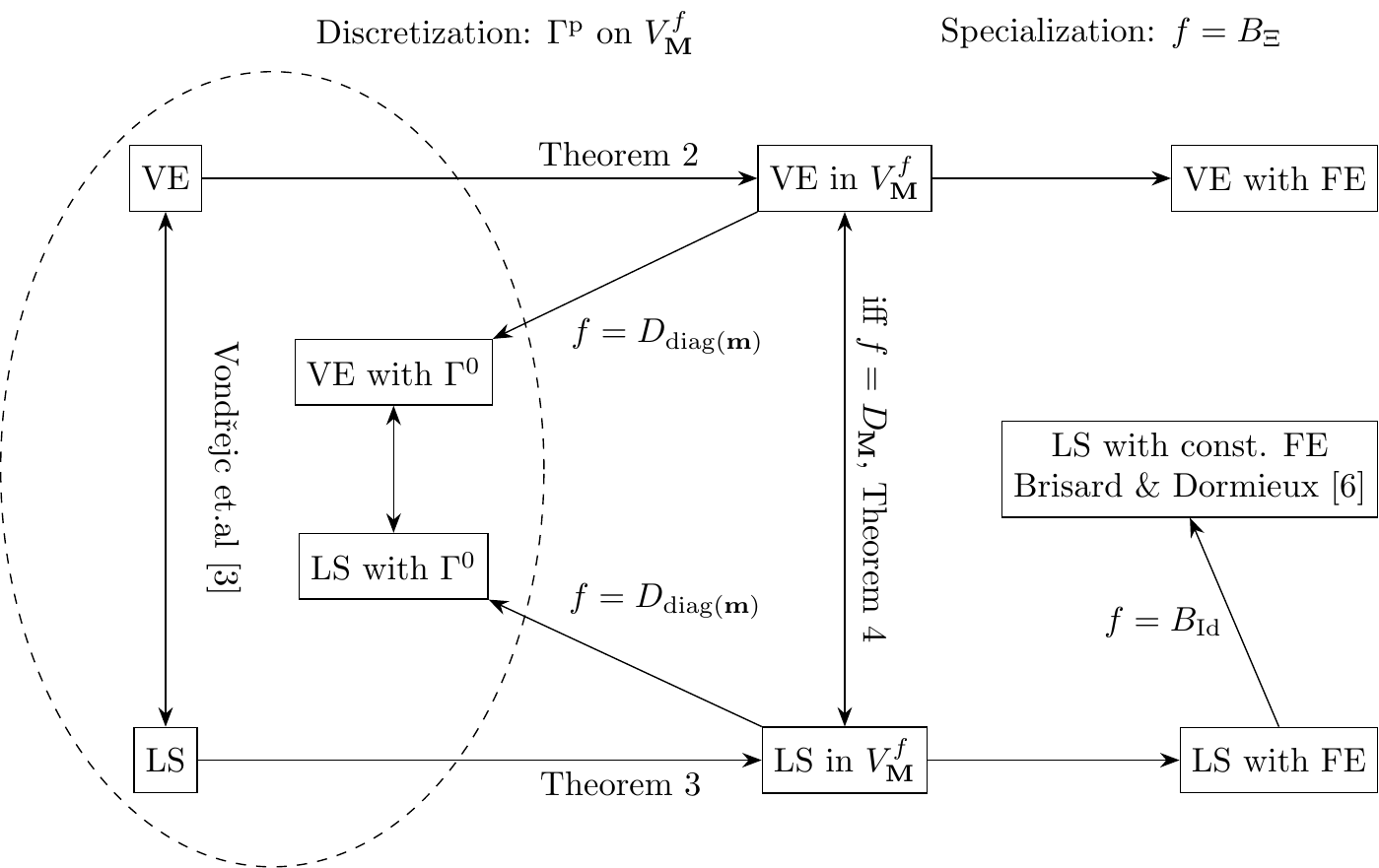}
  \caption{Illustration of the connections between solution spaces in this framework.}
  \label{fig:diagram}
\end{figure*}

In the following the discretisation of the variational
equation~\eqref{eq:weak_pde_projected} and the Lippmann-Schwinger
equation~\eqref{eq:lippmann-schwinger} on anisotropic translation invariant
spaces is introduced. The choice of the function $f$ generating the space
$V_\mat{M}^f$ influences properties of the discrete solution spaces. The
connections between these spaces are detailed in Figure~\ref{fig:diagram}. The
rest of this section follows this diagram and explains the relations therein.
Further details and proofs can be found in~\cite{BergmannMerkert2017}.

In the following, the discretised version of the strain $\varepsilon$ is denoted
by $\varepsilon_\mat{M} \in \bigl(V_\mat{M}^f\bigr)^{d \times d}$ with
$\varepsilon_\mat{M} = \sum_{\vect{y} \in \Pattern(\mat{M})} E_\vect{y}
\Translate(\vect{y}) \Fundamental_\mat{M}$. Test functions $\gamma_\mat{M}$ are also
expanded in terms of translates of $f$ by $\gamma_\mat{M} = \sum_{\vect{y} \in
\Pattern(\mat{M})} G_\vect{y} \Translate(\vect{y}) \Fundamental_\mat{M}$.
We further denote the discrete Fourier transforms of $\vect{E} =
(E_\vect{y})_{\vect{y} \in \Pattern(\mat{M})}$ and $\vect{G} =
(G_\vect{y})_{\vect{y} \in \Pattern(\mat{M})}$ by $\hat{\vect{E}} =
(\hat{E}_\vect{h})_{\vect{h} \in \generatingSet(\mat{M}^\tT)} =
\Fourier(\mat{M})\vect{E}$ and $\hat{\vect{G}} =
\Fourier(\mat{M})\vect{G}$, respectively.

The approximation of the solution of~\eqref{eq:weak_pde_projected}
and~\eqref{eq:lippmann-schwinger} requires the definition of the periodised
Green operator.

\begin{definition}
  Let $f \in L^2(\T^d)$, then the periodised Green operator $\Gamma^\p$ on
  $V_\mat{M}^f$  is defined in terms of its action on a symmetric field
  $\gamma_\mat{M} \in \bigl(V_\mat{M}^f\bigr)^{d \times d}$ by
  \begin{align}
    \label{eq:gamma_tilde_fourier}
    \Gamma^\p \tensorProd \gamma_\mat{M}
    &\coloneqq
    \sum_{\vect{y} \in \Pattern(\mat{M})} 
    \Gamma^\p_\vect{y} \tensorProd G_\vect{y} \Translate(\vect{y}) f.
    \intertext{In terms of Fourier sums this is the same as}
    \Gamma^\p \tensorProd \gamma_\mat{M}
    &\coloneqq
    \sum_{\vect{h} \in \generatingSet(\mat{M}^\tT)}
    \hat{\Gamma}^\p_\vect{h} \tensorProd \hat{G}_\vect{h}
    c_\vect{h}^\mat{M}(f) e^{2 \pi \im \vect{h}^\tT \cdot}
    \intertext{with Fourier coefficients}
    \hat{\Gamma}^\p_\vect{h} \tensorProd \hat{G}_\vect{h}
    &\coloneqq
    m \Bigl[
      \bigl\lbrace
      \hat{\Gamma}^0_\vect{k} \abs{c_\vect{k}(f)}^2
      \bigr\rbrace_{\vect{k} \in \Z^d}
    \Bigr]^\mat{M}_\vect{h} \tensorProd
    \hat{G}_\vect{h},\ \vect{h} \in \generatingSet(\mat{M}^\tT).
    \label{eq:bracket}
  \end{align}
\end{definition}

With these definitions the discretisation of the PDEs on spaces of translates
can be performed as follows:
\begin{theorem}
  Let the translates of $f$ be orthonormal and let $\Stiffness \in A(\T^d)^{d
  \times d \times d \times d}$ with the usual symmetries.
   Then $\varepsilon_\mat{M}$ fulfils the
  weak form
  \begin{equation}
    \label{eq:weak_pde_projected_1}
    \bigl\langle
    \varepsilon_\mat{M} +
    \Gamma^0 \bigl(
      \Stiffness - \Stiffness^0
    \bigr) \tensorProd
    \bigl(
      \varepsilon_\mat{M} + \varepsilon^0
    \bigr),
    \gamma_\mat{M}
    \bigr\rangle = 0
  \end{equation}
  for all symmetric $\gamma \in \bigl( V_\mat{M}^f \bigr)^{d \times d}$ if and only
  if the equation (VE in $V_\mat{M}^f$)
  \begin{equation}
    \label{eq:pde_coefficients}
    \sum_{\vect{y} \in \Pattern(\mat{M})}
    \Bigl(
      E_\vect{y} + 
      \Gamma^\p_\vect{y}
      \bigl(\Stiffness(\vect{y}) - \Stiffness^0 \bigr) \tensorProd
      \bigl( E_\vect{y} + \varepsilon^0 \bigr)
    \Bigr)
    \bigl(\Translate(\vect{y}) \Fundamental_\mat{M} \bigr)(\vect{x}) = \vect{0}
  \end{equation}
  holds true for all $\vect{x} \in \Pattern(\mat{M})$.
\end{theorem}

\begin{theorem}
  Let the translates of $f$ be orthonormal and let $\Stiffness \in A(\T^d)^{d
  \times d \times d \times d}$ with the usual symmetries.
   Then $\varepsilon_\mat{M}$ fulfils the
  weak form
  \begin{equation}
    \label{eq:weak_pde_projected_2}
    \bigl\langle
    \Stiffness^0
    \Gamma^0
      \Stiffness
    \tensorProd
    \bigl(
      \varepsilon_\mat{M} + \varepsilon^0
    \bigr),
    \gamma_\mat{M}
    \bigr\rangle = 0
  \end{equation}
  for all symmetric $\gamma \in \bigl( V_\mat{M}^f \bigr)^{d \times d}$ if and only
  if the equation (LS in $V_\mat{M}^f$)
  \begin{equation}
    \label{eq:pde_coefficients_variational}
    \sum_{\vect{y} \in \Pattern(\mat{M})}
      \Stiffness^0 \Gamma^\p_\vect{y}
      \Stiffness(\vect{y}) \tensorProd
      \bigl( E_\vect{y} + \varepsilon^0 \bigr)
      \bigl(\Translate(\vect{y}) \Fundamental_\mat{M} \bigr)(\vect{x}) = \vect{0}
  \end{equation}
  holds true for all $\vect{x} \in \Pattern(\mat{M})$.
\end{theorem}

Choosing the function $f$ as the Dirichlet kernel $f = D_\mat{M}$
reduces~\eqref{eq:pde_coefficients} and~\eqref{eq:pde_coefficients_variational}
to discretised equations which correspond to truncating the Fourier series
arising in the continuous equations~\eqref{eq:weak_pde_projected}
and~\eqref{eq:lippmann-schwinger}. The effect of the pattern matrix $\mat{M}$ in
this special case is analysed in~\cite{BergmannMerkert2016}. For the Dirichlet
kernel on a tensor product grid, which corresponds to a diagonal pattern matrix,
Vond\v{r}ejc et.al.~\cite{Vondrejc2014} prove the equivalence of the continuous
equations and show that it carries over to the discrete solutions as well. The
equations (VE) and (LS) discretised on such a tensor product grid using the
Dirichlet kernel, i.e.~$f = D_{\diag(\vect{m})}$ for $\vect{m} \in \N^d$ are
denoted by (VE with $\Gamma^0$) and (LS with $\Gamma^0$) in the diagram,
respectively.

The Dirichlet kernel proves to play a special role when analysing the
equivalence of the discrete solutions.

\begin{theorem}
  Let the translates of $f$ be orthonormal and let $\Stiffness \in A(\T^d)^{d
  \times d \times d \times d}$ with the usual symmetries.
  Then $\varepsilon_\mat{M}$ fulfils the equations
  \begin{equation}
    \sum_{\vect{y} \in \Pattern(\mat{M})}
    \Bigl(
      E_\vect{y} + 
      \Gamma^\p_\vect{y}
      \bigl(\Stiffness(\vect{y}) - \Stiffness^0 \bigr) \tensorProd
      \bigl( E_\vect{y} + \varepsilon^0 \bigr)
    \Bigr)
    \bigl(\Translate(\vect{y}) \Fundamental_\mat{M} \bigr)(\vect{x}) = \vect{0}
  \end{equation}
  and
  \begin{equation}
    \sum_{\vect{y} \in \Pattern(\mat{M})}
      \Stiffness^0 \Gamma^\p_\vect{y}
      \Stiffness(\vect{y}) \tensorProd
      \bigl( E_\vect{y} + \varepsilon^0 \bigr)
      \bigl(\Translate(\vect{y}) \Fundamental_\mat{M} \bigr)(\vect{x}) = \vect{0}
  \end{equation}
  for all $\vect{x} \in \Pattern(\mat{M})$ if and only if $V_\mat{M}^f =
  V_\mat{M}^{D_\mat{M}}$.
\end{theorem}

The equivalence of the solutions in case of the space generated by the
Dirichlet kernel is due to the representation of constant functions in frequency
domain. In case of the Dirichlet kernel, i.e. in case of
the discrete Fourier transform, a constant
function $g$ is only characterized by the Fourier coefficient $c_\vect{0}(g)$.
This is no longer true for other spaces of translates which leads to different
solutions of~\eqref{eq:pde_coefficients} and~\eqref{eq:weak_pde_projected_2}.

Functions $f$ with compact support in frequency domain allow for an easy
evaluation of the Bracket sum in~\eqref{eq:bracket}. This is no longer the case
when choosing $f$ compactly supported in space domain. However, such
functions allow to include (simplified) finite elements into the framework by
choosing $f$ to be a Box spline on the unit cell $\mat{M}^{-1}
\bigl[-\frac12,\frac12\bigr)^d$. The case of constant Box splines was already
analysed by Brisard and Dormieux~\cite{BrisardDormieux2010Framework} in case of
the Lippmann-Schwinger equation.

\begin{figure*}[t]
  \centering
  \begin{subfigure}[c]{.32\textwidth}\centering
    \includegraphics[width=.98\textwidth]{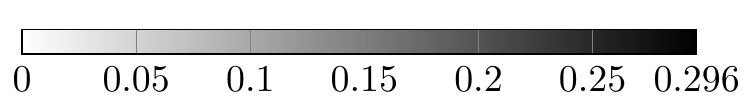}
  \end{subfigure}\\[.25\baselineskip]
  \begin{subfigure}[c]{.32\textwidth}\centering
    \includegraphics[width=.98\textwidth]{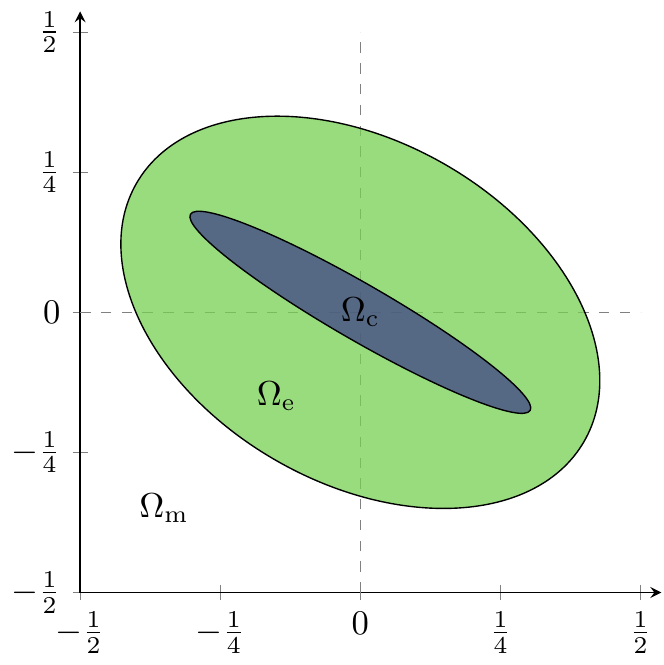}
  \end{subfigure}%
  \begin{subfigure}[c]{.32\textwidth}\centering
    \includegraphics[width=.98\textwidth]{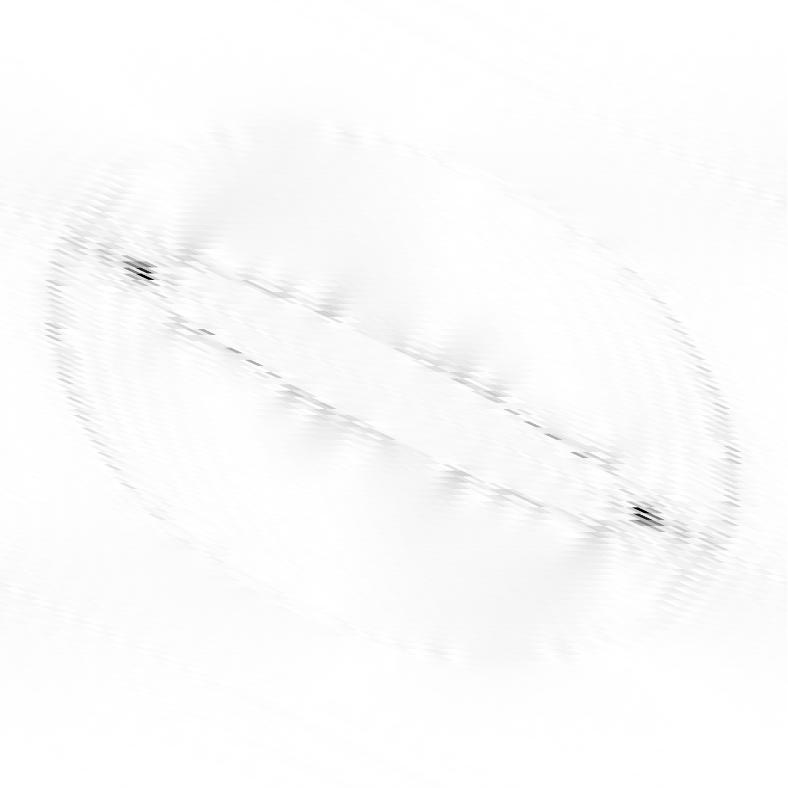}\\[1mm]
    \small
    \begin{tabular}{c l}
      \multirow{2}{*}{$f=D_{\mat{M}}$,} & $e_\eff = 0.0036$,\\
      &$e_{\ell^2} = 0.022$.
    \end{tabular}
  \end{subfigure}%
  \begin{subfigure}[c]{.32\textwidth}\centering
    \includegraphics[width=.98\textwidth]{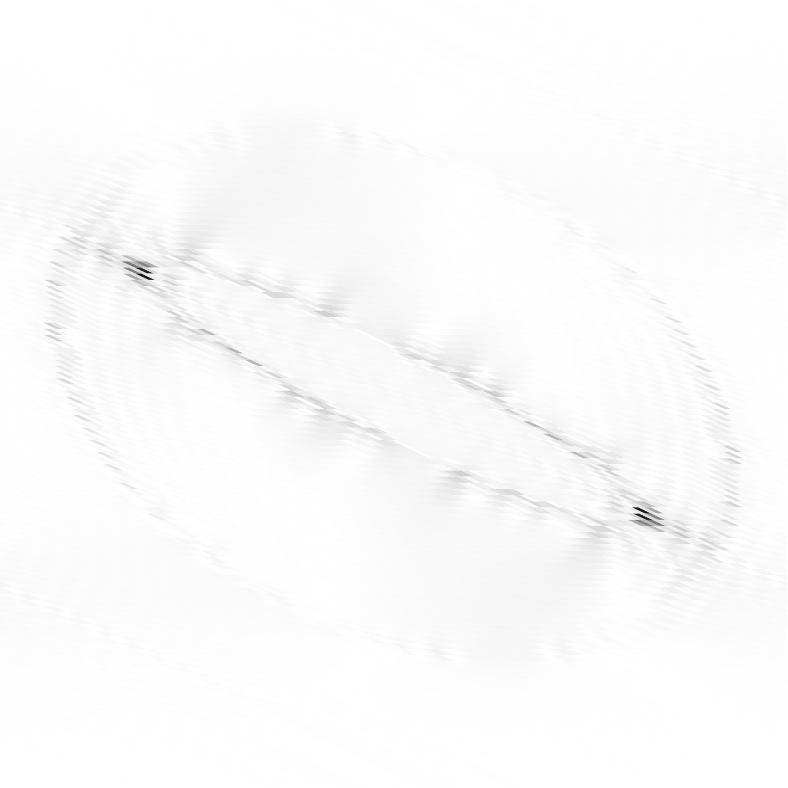}\\[1mm]
    \small
    \begin{tabular}{c l}
      \multirow{3}{*}{$f=f_{\mat{M},\alpha}$,} & $\alpha = (0.4, 0)^\tT$,\\
      & $e_\eff = 0.0024$,\\
      &$e_{\ell^2} = 0.025$.
    \end{tabular}
  \end{subfigure}
  \caption{A schematic of the generalized Hashin structure (left) and the
  $e_{\log}$ error using the Dirichlet kernel (middle), and a de la Vall\'ee
Poussin mean (right) based solution.}
  \label{fig:hashin}
\end{figure*}

\section{Numerical examples}\label{sec:numerics}

The effect of choosing different spaces of translates $V_\mat{M}^f$ is
demonstrated in the following. The generalized Hashin structure that consists of
two confocal ellipsoids ($\Omega_c$) and ($\Omega_e$) embedded in a matrix
material ($\Omega_m$) and is depicted in Figure~\ref{fig:hashin} (left). For
this structure an analytic expression for the strain $\varepsilon$ and for the
effective stiffness $\Stiffness^\eff \tensorProd \varepsilon^0$ is known and
described in~\cite[Section 4.2]{BergmannMerkert2016}. The examples were computed
using~\eqref{eq:pde_coefficients} using an iterative scheme based on a Neumann
series approach, originally proposed by Moulinec and
Suquet~\cite{MoulinecSuquet1998}.

The following computation was performed with the pattern matrix $\mat{M} =
\bigl( \begin{smallmatrix} 128 & 272\\ 0 & 128 \end{smallmatrix}\bigr)$, showing
the minimal $\ell^2$-error in~\cite{BergmannMerkert2016}.

The relative logarithmic error $e_{\log} = \log(1+\abs{\varepsilon +
\tilde\varepsilon})$ with $\tilde\varepsilon$ being the analytic solution is
depicted for the Dirichlet kernel (middle) and de la Vall\'ee Poussin means with
slopes~$\alpha = (0.4,0)^\tT$ (right). The error using de la Vall\'ee Poussin
means shows significantly less Gibbs phenomenon and thus a smoother solution.
The benefits of choosing~$f=f_{\mat{M},\alpha}$ also get visible when looking at
the $\ell^2$-error and the error in the effective stiffness.
These errors are defined via $e_{\ell^2} = \norm{\varepsilon-\tilde\varepsilon}
\norm{\tilde\varepsilon}^{-1}$ and $e_\eff = \norm{ \tilde\Stiffness^\eff
\tensorProd \varepsilon^0 - \Stiffness^\eff \tensorProd \varepsilon^0}
\norm{\tilde\Stiffness^\eff \tensorProd \varepsilon^0}^{-1}$, respectively. The
analytic effective stiffness matrix is denoted by $\tilde\Stiffness^\eff$.

The slopes $\alpha$ of the de la Vall\'ee Poussin means were optimized to yield
the smallest $e_\eff$ which is reduced by more than $33\%$. The $\ell^2$-error
increases slightly from $0.022$ to $0.025$.

\section{Conclusion}\label{sec:conclusion}

The introduced framework extends the truncated Fourier series approach to
anisotropic spaces of translates and thus incorporates amongst others simplified
finite elements. The Green operator $\Gamma^0$ emerges as a special case of
Dirichlet kernel translates. Especially, the discrete solutions of the
Lippmann-Schwinger equation and the variational equation coincide if and only if
the space of translates $V_\mat{M}^f$ is the one generated by the Dirichlet
kernel. Choosing de la Vall\'ee Poussin means for $f$ leads to a reductions of
the Gibbs phenomenon and allows for a better prediction of the effective
stiffness of a composite material.

Continuing the work of~\cite{BergmannMerkert2017} a direct extension towards
periodic wavelets and finite elements with full integration are two points of
future work. Furthermore, the convergence of the discretisation towards the
continuous solution and convergence speed have to be investigated.

\printbibliography
\end{document}